\documentclass{IEEEtran4PSCC}
\ifCLASSINFOpdf
   \usepackage[pdftex]{graphicx}
\else
   \usepackage[dvips]{graphicx}
\fi
\usepackage[cmex10]{amsmath}
\usepackage{textcomp}
\usepackage{subfig}
\usepackage{lineno}
\usepackage[hidelinks]{hyperref}
\usepackage{mathtools}
\usepackage{amsthm}

\usepackage[dvipsnames]{xcolor}
\definecolor{dkred}{rgb}{0.8,0,0}
\definecolor{blue}{rgb}{0,0,1}

\usepackage{pgfplots}
\pgfplotsset{compat=1.18}
\usepackage{multirow}
\usepackage{booktabs}
\usepackage{tabularx}
\usepackage{xcolor}
\usepackage{tikz}
\usetikzlibrary{patterns}
\usepackage{amssymb}
\usepackage{graphicx}
\usepackage{wasysym}
\usepackage{enumitem}
\usepackage[table]{xcolor}
\usepackage{pgfplots}
\pgfplotsset{compat=1.18}
\usetikzlibrary{positioning}
\usetikzlibrary{spy}
\usepgfplotslibrary{statistics, groupplots}
\usepackage{mathrsfs}
\usepackage{makecell}
\usepackage{graphicx}
\usepackage{textcomp}
\usepackage{algorithm}
\usepackage{algpseudocode}  
\usepackage{float}
\usetikzlibrary{arrows.meta}
\usepackage{nicematrix}
\usepackage{arydshln}
\usepackage{cuted}
\usepackage{adjustbox}
\usepackage{comment}

\makeatletter
\let\old@ps@headings\ps@headings
\let\old@ps@IEEEtitlepagestyle\ps@IEEEtitlepagestyle
\def\psccfooter#1{%
    \def\ps@headings{%
        \old@ps@headings%
        \def\@oddfoot{\strut\hfill#1\hfill\strut}%
        \def\@evenfoot{\strut\hfill#1\hfill\strut}%
    }%
    \def\ps@IEEEtitlepagestyle{%
        \old@ps@IEEEtitlepagestyle%
        \def\@oddfoot{\strut\hfill#1\hfill\strut}%
        \def\@evenfoot{\strut\hfill#1\hfill\strut}%
    }%
    \ps@headings%
}
\makeatother

\begin{document}

\title{Towards Exact Temporal Aggregation of Time-Coupled Energy Storage Models via Active Constraint Set Identification and Machine Learning}

\author{
\IEEEauthorblockN{Thomas Klatzer, David Cardona-Vasquez, Luca Santosuosso and Sonja Wogrin}
\IEEEauthorblockA{Institute of Electricity Economics and Energy Innovation, \\ 
Research Center for Energy Economics and Energy Analytics \\
Graz University of Technology\\
Graz, Austria\\
\{thomas.klatzer, david.cardonavasquez, luca.santosuosso, wogrin\}@tugraz.at}
}

\maketitle

\begin{abstract}
Time series aggregation (TSA) aims to construct temporally aggregated optimization models that accurately represent the output space of their full-scale counterparts while using a significantly reduced temporal dimensionality. This paper presents a theoretical approach that achieves exact temporal aggregation of full-scale power system models -- even in the presence of energy storage time-coupling constraints -- by leveraging active constraint sets and dual information. This advances the state of the art beyond existing TSA methods, which typically cannot guarantee solution accuracy or rely on iterative procedures to determine the required number of representative periods. To bridge the gap between this theoretical analysis and practical application, we employ machine learning, i.e., classification and clustering, to inform TSA in models that co-schedule variable renewable energy sources and energy storage. Numerical results show substantially improved computational performance relative to the full-scale model, while maintaining a favorable trade-off between solution accuracy and complexity.
\end{abstract}

\begin{IEEEkeywords}
Computational complexity, energy storage system, machine learning, optimization, time series aggregation. \relax
\end{IEEEkeywords}

\thanksto{Funded by the European Union (ERC, NetZero-Opt, 101116212). Views and opinions expressed are however those of the authors only and do not necessarily reflect those of the European Union or the European Research Council. Neither the European Union nor the granting authority can be held responsible for them.}

\section*{Nomenclature}
\subsection{Sets}
\begin{IEEEdescription}[\IEEEusemathlabelsep\IEEEsetlabelwidth{$\smash{\overline{P}}^\mathrm{c}_s, \smash{\overline{P}}^\mathrm{d}_s$}]
\item[$\mathcal{G}$]                Set of all generators $g$.
\item[$\mathcal{G}^\mathcal{V}$]    Subset of renewable generators $v$.
\item[$\mathcal{G}^\mathcal{T}$]    Subset of thermal generators $t$.
\item[$\mathcal{S}$]                Set of energy storage systems $s$.
\item[$\mathcal{I}$]                Set of independent submodels $i$.
\item[$\mathcal{R^\mathcal{I}}$]    Set of periods $r$ of independent submodel $i$.
\end{IEEEdescription}

\subsection{Parameters}
\begin{IEEEdescription}[\IEEEusemathlabelsep\IEEEsetlabelwidth{$\smash{\overline{P}}^\mathrm{c}_s, \smash{\overline{P}}^\mathrm{d}_s$}]
\item[$W_r$]                        Weight of period $r$ (h).
\item[$C_g$]                        Variable operating cost of $g$ (\texteuro/MWh).
\item[$C^\mathrm{d}_s$]             Variable cost of discharging $s$ (\texteuro/MWh).
\item[$C^\mathrm{ns}$]              Cost of non-supplied energy (\texteuro/MWh).
\item[$\widetilde{D}_r$]            Average power demand in $r$ (MW).
\item[$\widetilde{F}_{r,v}$]        Average capacity factor of $v$ in $r$ (p.u.).
\item[$\overline{P}_g$]             Max. power generation of $g$ (MW).
\item[$\smash{\overline{P}}^\mathrm{c}_s, \smash{\overline{P}}^\mathrm{d}_s$]   Max. charging/discharging power of $s$ (MW).
\item[$\eta^\mathrm{c}_s, \eta^\mathrm{d}_s$]                                   Charging/discharging efficiency of $s$ (--).
\item[$\underline{E}_s,\overline{E}_s$]                                         Min./max. state of charge of $s$ (MWh).
\end{IEEEdescription}

\subsection{Primal Variables}
\begin{IEEEdescription}[\IEEEusemathlabelsep\IEEEsetlabelwidth{$\smash{\overline{P}}^\mathrm{c}_s, \smash{\overline{P}}^\mathrm{d}_s$}]
\item[$p_{r,g}$]                                Power generation of $g$ (MW).
\item[$p^\mathrm{ns}_{r}$]                      Non-supplied power (MW).
\item[$p^\mathrm{c}_{r,s},p^\mathrm{d}_{r,s}$]  Charging/discharging power of $s$ (MW).
\item[$e_{r,s}$]                                State of charge of $s$ (MWh).
\end{IEEEdescription}

\subsection{Dual Variables}
\begin{IEEEdescription}[\IEEEusemathlabelsep\IEEEsetlabelwidth{$\smash{\lambda}^\mathrm{ini}_{r,s},\smash{\lambda}^\mathrm{soc}_{r,s},\smash{\lambda}^\mathrm{end}_{r,s}$}]
\item[$\smash{\mu}^\mathrm{bal}_r$]                                                                      Dual of the power balance constraint.
\item[$\smash{\underline{\lambda}}^\mathrm{v}_{r,v}, \smash{\overline{\lambda}}^\mathrm{v}_{r,v}$]       Duals of lower/upper bounds of $p_{r,v}$.
\item[$\smash{\underline{\lambda}}^\mathrm{t}_{r,t}, \smash{\overline{\lambda}}^\mathrm{t}_{r,t}$]       Duals of lower/upper bounds of $p_{r,t}$.
\item[$\smash{\underline{\lambda}}^\mathrm{ns}_{r}, \smash{\overline{\lambda}}^\mathrm{ns}_{r}$]         Duals of lower/upper bounds of $p^\mathrm{ns}_r$.
\item[$\smash{\mu}^\mathrm{ini}_{r,s}, \smash{\mu}^\mathrm{intra}_{r,s}, \smash{\mu}^\mathrm{fin}_{r,s}$] Duals of the state-of-charge constraints.
\item[$\smash{\underline{\lambda}}^\mathrm{soc}_{r,s}, \smash{\overline{\lambda}}^\mathrm{soc}_{r,s}$]   Duals of lower/upper bounds of $e_{r,s}$.
\item[$\smash{\underline{\lambda}}^\mathrm{c}_{r,s}, \smash{\overline{\lambda}}^\mathrm{c}_{r,s}$]       Duals of lower/upper bounds of $p^\mathrm{c}_{r,s}$.
\item[$\smash{\underline{\lambda}}^\mathrm{d}_{r,s}, \smash{\overline{\lambda}}^\mathrm{d}_{r,s}$]       Duals of lower/upper bounds of $p^\mathrm{d}_{r,s}$.

\end{IEEEdescription}

\section{Introduction}
\label{sec:Intro}
Energy storage systems (ESSs), such as batteries and hydro storage, are critical for mitigating the stochasticity of variable renewable energy (VRE) sources~\cite{DeCarne2024}. However, co-scheduling VRE sources and ESSs over extended time horizons or at high temporal resolution presents significant computational challenges, arising from the high-dimensional solution space of the optimization model, the time-coupling constraints of ESSs, and the intricate interactions among heterogeneous units~\cite{santosuosso2024distributed}.
This complexity drives the adoption of time series aggregation (TSA) methods, which substantially reduce input time series dimensionality, yielding an aggregated model that aims to accurately approximate its full-scale counterpart~\cite{Hoffmann2020}.

Standard a priori TSA methods aim to capture the statistical features of the input time series, such as demand and VRE generation, while reducing dimensionality through clustering techniques including $k$\nobreakdash-means~\cite{Almaimouni2018,Liu2018,Kotzur2018a}, $k$\nobreakdash-medoids~\cite{Kotzur2018a,Tejada-Arango2018,Kotzur2018}, or hierarchical clustering~\cite{Liu2018,Kotzur2018a,Nahmmacher2016,Pineda2018}. However, prior work~\cite{Teichgraeber2022} has shown that accurately capturing the input space of a model does not necessarily translate into accurately capturing its output space, i.e., its solution space. Instead, accounting for the system's structural characteristics during TSA, particularly ESS time-coupling constraints, is crucial for accurately approximating the full-scale model~\cite{Kotzur2018a}. 

To overcome the limitations of traditional a priori TSA, the concept of a posteriori TSA has recently emerged, focused on minimizing errors in the aggregated model's output space by incorporating model information into the TSA~\cite{Hoffmann2020}.
For instance, the TSA method in Sun et al.~\cite{Sun2019} relies on the objective function values of daily partitions of the full time horizon, which are clustered using hierarchical clustering. Li et al.~\cite{Li2022a} follow a similar approach, employing $k$-means and $k$-medoids clustering. In contrast, Teichgraeber et al.~\cite{Teichgraeber2020} iteratively add days for which no feasible model solution exists to the set of representative periods. The TSA method in Bahl et al.~\cite{Bahl2017} starts from an initial clustering-based set of representative periods, gradually adds infeasible periods, and re-performs clustering if the objective function-based termination criterion cannot be satisfied by the aggregated model. In another work, Teichgraeber and Brandt~\cite{Teichgraeber2021} classify days with the highest lost load as extreme and add them iteratively until lost load is reduced to an acceptable threshold. Hilbers et al.~\cite{Hilbers2023} similarly identify days with high generation costs as extreme and further incorporate storage charging and discharging patterns from previous solutions into TSA.

From the literature, it is evident that existing a posteriori TSA methods typically rely on expert-defined, case-study-dependent partitions of the time horizon (e.g., into representative days).
However, these methods typically do not provide theoretical guarantees on the quality of the aggregated model’s output relative to its full-scale counterpart and thus remain largely heuristic.
When performance-guaranteed TSA methods are proposed -- typically providing bounds on the maximum approximation error incurred by the aggregated model~\cite{santosuosso2025optimal} -- they require iterative procedures to converge on the number of representative periods needed to adequately capture full-scale model dynamics (e.g.,~\cite{Bahl2017}).
In contrast, Wogrin~\cite{Wogrin2023} proposed an a posteriori TSA method that leverages the identification of active constraints in the full-scale model to perform TSA and construct an aggregated model that \textit{exactly} preserves the optimal solution space. Yet, this method does not support models with ESS time-coupling constraints and relies on ex-ante knowledge of the full-scale model solution for TSA.

Against this background, the main contributions of this paper are as follows: 
(1) We extend the theoretical results from~\cite{Wogrin2023}, demonstrating that \textit{exact} TSA is theoretically achievable in models involving ESS time-coupling constraints. This is accomplished through the identification of active constraint sets (ACSs)\footnote{A period's ACS consists of all constraints that are binding (satisfied with equality) at a given solution point, thereby defining the boundary of the feasible region where an optimal solution lies.}, which inform the disaggregation of the full-scale model into parallelizable submodels, and the subsequent aggregation of consecutive periods within each submodel using dual information. 
(2) Given that ACSs and corresponding dual information stems from ex-ante knowledge of the full-scale model solution, we propose a machine learning (ML) framework, i.e., classification and clustering, to inform TSA, thereby bridging the provided theoretical result with practical applicability of the novel TSA method.
(3) We validate the ML framework across different VRE and ESS configurations -- such as solar or wind paired with battery or pumped hydro storage -- to capture their distinct operational dynamics and technical characteristics, and analyze their impacts on parallelization and aggregation potential.

The remainder of the paper is organized as follows: Section~\ref{sec:Model} presents the co-scheduling model allowing for exact TSA. Section~\ref{sec:ExactTSA} demonstrates that exact TSA is theoretically achievable when the optimization model includes ESS time-coupling constraints. Building on this theoretical results, Section~\ref{sec:ML} introduces the ML framework for disaggregation and subsequent aggregation, which is validated through numerical experiments in Section~\ref{sec:NumRes}. Finally, Section~\ref{sec:Conclusions} concludes the paper and outlines directions for future research.

\section{Optimization Model for Exact TSA with ESSs}
\label{sec:Model}
\allowdisplaybreaks

This section introduces the optimization model~\eqref{eqn:model}, which features a flexible temporal structure that enables TSA through disaggregation and subsequent aggregation, thereby exactly preserving the optimal solution of its full-scale counterpart.
This approach crucially enables parallelization of the constructed aggregated models, which reduces computational complexity and improves scalability. We achieve this by disaggregating the full time horizon into independent submodels $i\in\mathcal{I}$, each defined over its own set of representative periods $\mathcal{R}^i$, indexed by $r$. For example, a time horizon of 8736 periods can be represented using a single submodel $|\mathcal{I}|=1$ with $|\mathcal{R}^1|=8736$ periods\footnote{We use $|\cdot|$ to denote the cardinality of a set.}, or by any other meaningful combination of submodels and periods.
In addition to disaggregating the full time horizon into independent submodels, periods within each submodel can be aggregated. To this end, each $r\in\mathcal{R}^i$ is assigned a weight $W_{r}$, representing the number of periods aggregated into~$r$. For example, a submodel representing a single day may consist of 24 periods with weights $W_{r}=1$, or a single representative period with $W_{r}=24$, or any other combination of weights and periods. For each representative period $r$, the input time series are represented by the average over the corresponding time interval. 

For each submodel~$i$, the objective function~\eqref{eqn:of} minimizes the total operating costs of VRE and thermal units, ESSs, and non-supplied power, with associated costs $C_g, C^{\mathrm{d}}_s$ and $C^{\mathrm{ns}}$ respectively. By summing the individual objective function values~$z_i$ of each submodel, we derive the total objective incurred over the full time horizon of the co-scheduling problem. Constraint~\eqref{eqn:balance} presents the power balance. Generation from thermal units $p_{r,t}$ is limited by their installed capacity $\overline{P}_t$, as shown in~\eqref{eqn:therm_bounds}. Generation from VRE sources $p_{r,v}$ is constrained by their installed capacity~$\overline{P}_v$ and the average capacity factor~$\widetilde{F}_{r,v}$, which represents the available potential in each period, as given in~\eqref{eqn:res_bounds}. Finally, non-supplied power~$p^{\mathrm{ns}}_{r}$ is limited by the average demand~$\widetilde{D}_{r}$, as defined in~\eqref{eqn:nsp_bounds}.
Constraints~\eqref{eqn:soc_start}--\eqref{eqn:soc_fin} govern the evolution of the storage state of charge~$e_{r,s}$ across the initial~\eqref{eqn:soc_start}, intra-period~\eqref{eqn:soc_intra}, and final~\eqref{eqn:soc_fin} periods of the submodel, with charging~$p^\mathrm{c}_{r,s}$ and discharging~$p^\mathrm{d}_{r,s}$ accounted for through their respective efficiencies, $\eta^\mathrm{c}_{s}$ and $\eta^\mathrm{d}_{s}$, and weighted by $W_r$ to convert power to energy.
Finally, the state of charge is constrained by its lower and upper limits,~$\underline{E}_s$ and~$\overline{E}_s$, as given in~\eqref{eqn:soc_bounds}, while~$p^\mathrm{c}_{r,s}$ and~$p^\mathrm{d}_{r,s}$ are limited by their respective capacities~$\smash{\overline{P}}^\mathrm{c}_s$ and~$\smash{\overline{P}}^\mathrm{d}_s$ in~\eqref{eqn:sto_ch} and~\eqref{eqn:sto_dis}. Since the co-scheduling problem is formulated as a linear program, dual variables exist and are denoted by Greek letters in parentheses.
\begin{figure*}[t]
\centering
\resizebox{\textwidth}{!}{\definecolor{Bars}{HTML}{FFC44D}

\begin{tikzpicture}[node distance=2mm]

\node[
    draw=none,
    rounded corners=0.5mm,
    fill=gray!40,
    minimum width=19.2cm,
    minimum height=1.9cm,
    anchor=north
] (box1) at (0,0) {};

\node[anchor=north] at (box1.north) {%
    \begin{minipage}{19cm}
        Step 1: Compute full-scale model~\eqref{eqn:model} and identify \textbf{ACSs and corresponding dual information} for each period.
    \end{minipage}%
};

\node[
    draw=none,
    rounded corners=0.5mm,
    fill=gray!15,
    minimum width=19cm,
    minimum height=1.2cm,
    anchor=south
] (box1b) at ([yshift=+1mm]box1.south) {}; 

\node[anchor=north] at (box1b.north) {%
    \begin{minipage}{18.8cm}
        \begin{itemize}[itemsep=4pt, topsep=4pt, left=5pt, label=\raisebox{0.35ex}{\scalebox{0.4}{$\blacksquare$}}]
            \item \textbf{Full-scale model}: $|\mathcal{I}|=1,~|\mathcal{R}^1|=8736,~W_r = 1$
            \item \textbf{Computational burden}: $187 \times 10^{-3}$~work units
        \end{itemize}
    \end{minipage}%
};

\node[
    draw=none,
    rounded corners=0.5mm,
    fill=gray!40,
    below=1mm of box1.south,
    minimum width=19.2cm,
    minimum height=4.7cm,
    anchor=north
] (box2) {};

\node[anchor=north] at (box2.north) {%
    \begin{minipage}{19cm}
        Step 2: \textbf{Disaggregate} full-scale model into independent \textbf{unlinked} submodels $\mathcal{U}$ and \textbf{linked} submodels $\mathcal{L}$ such that $\mathcal{I}=\mathcal{U}\cup\mathcal{L}$.
    \end{minipage}%
};

\node[
    draw=none,
    rounded corners=0.5mm,
    fill=gray!15,
    minimum width=19cm,
    minimum height=4cm,
    anchor=south
] (box2b) at ([yshift=+1mm]box2.south) {};

\node[anchor=north] at (box2b.north) {%
    \begin{minipage}{18.8cm}
        \begin{itemize}[itemsep=4pt, topsep=4pt, left=5pt, label=\raisebox{0.35ex}{\scalebox{0.4}{$\blacksquare$}}]
            \item \textbf{Unlinked} submodels represent 1377 periods using: $i\in\mathcal{U}=\{1,2,\dots,1377\},~|\mathcal{R}^i|\in\mathcal{R}^\mathcal{U}=\{1,1,\dots,1\},~W_r=1$
            \item \textbf{Linked} submodels represent 7359 periods using: $i\in\mathcal{L}=\{1,2\dots,176\},~|\mathcal{R}^i|\in\mathcal{R}^\mathcal{L}=\{4,5,\dots,231,253,395\},~W_r=1$
            \item \textbf{Lower bound on computational burden}: $77 \times 10^{-4}$~work units~$\rightarrow$~\textbf{24-fold} reduction versus full-scale model
        \end{itemize}
    \end{minipage}%
};

\node[anchor=center, inner sep=0,yshift=-0.9cm] at (box2b.center) {%
\begin{tikzpicture}
\begin{groupplot}[
    group style={
        group size=2 by 1,
        horizontal sep=1cm,
    },
    scale only axis,
    height=1.3cm,
    ybar interval,
    ymin=0,
    xtick=data,
    xticklabel={%
        \pgfmathtruncatemacro{\a}{\tick}%
        \pgfmathtruncatemacro{\b}{\nexttick-1}%
        \a--\b%
    },
    yticklabel style={font=\tiny, xshift=-2pt},
    xticklabel style={yshift=1.1cm, font=\tiny, anchor=east, rotate=90},
    tick style={major tick length=2pt},
    xtick style={draw=none},
    grid=none,
    axis on top,
]

\nextgroupplot[
    width=0.9cm,
    font=\scriptsize,
    title={Unlinked},
    title style={font=\scriptsize, yshift=-2pt},
    ylabel={Submodels},
    xlabel={$|\mathcal{R}^i|\in\mathcal{R}^\mathcal{U}$},
    xlabel style={font=\scriptsize, yshift=-2pt},
    ylabel style={font=\scriptsize, yshift=5pt},
    xmax=2,
    enlarge x limits=0.9,
]
\addplot+[
    hist={data=y, bins=1},
    fill=Bars,
    draw=black,
    mark=none
] table[y index=0, col sep=comma] {method/histLength_unlinked.csv};

\nextgroupplot[
    width=10cm,
    font=\scriptsize,
    title={Linked},
    title style={font=\scriptsize, yshift=-2pt},
    xlabel={$|\mathcal{R}^i|\in\mathcal{R}^\mathcal{L}$},
    xlabel style={font=\scriptsize, yshift=-2pt},
    ymax=90,
    enlarge x limits=0.02,
]
\addplot+[
    hist={data=y, bins=30, data min=4, data max=396},
    fill=Bars,
    draw=black,
    mark=none
] table[y index=0, col sep=comma] {method/histLength_linked.csv};

\end{groupplot}
\end{tikzpicture}%
};

\node[
    draw=none,
    rounded corners=0.5mm,
    fill=gray!40,
    below=1mm of box2.south,
    minimum width=19.2cm,
    minimum height=1.9cm,
    anchor=north
] (box3) {};

\node[anchor=north] at (box3.north) {%
    \begin{minipage}{19cm}
        Step 3: \textbf{Aggregate} \textbf{unlinked} submodels, and aggregate adjacent periods within \textbf{linked} submodels sharing the same ACS and duals.
    \end{minipage}%
};

\node[
    draw=none,
    rounded corners=0.5mm,
    fill=gray!15,
    minimum width=19cm,
    minimum height=1.2cm,
    anchor=south
] (box3b) at ([yshift=+1mm]box3.south) {};

\node[anchor=north] at (box3b.north) {%
    \begin{minipage}{18.8cm}
        \begin{itemize}[itemsep=4pt, topsep=4pt, left=5pt, label=\raisebox{0.35ex}{\scalebox{0.4}{$\blacksquare$}}]
            \item \textbf{Unlinked} submodels aggregate into $|\mathcal{R}^i|\in\mathcal{R}^\mathcal{U}=\{1,1,\dots,1\}$ with $W_r\geq1~\rightarrow~\sum_{|\mathcal{R}^i|\in\mathcal{R}^\mathcal{U}}=8$ periods
            \item \textbf{Linked} submodels aggregate into $|\mathcal{R}^i|\in\mathcal{R}^\mathcal{L}=\{2,2,\dots,31,27,5\}$ with $W_r\geq1~\rightarrow~\sum_{|\mathcal{R}^i|\in\mathcal{R}^\mathcal{L}}=925$ periods 
        \end{itemize}
    \end{minipage}%
};

\node[
    draw=none,
    rounded corners=0.5mm,
    fill=gray!40,
    below=1mm of box3.south,
    minimum width=19.2cm,
    minimum height=1.8cm,
    anchor=north
] (box4) {};

\node[anchor=north] at (box4.north) {%
    \begin{minipage}{19cm}
        Step 4: \textbf{Solve disaggregated and subsequently aggregated submodels} in parallel and compute total OFV $Z = \sum_{i\in\mathcal{I}} z_i$.
    \end{minipage}%
};

\node[
    draw=none,
    rounded corners=0.5mm,
    fill=gray!15,
    minimum width=19cm,
    minimum height=1.1cm,
    anchor=south
] (box4b) at ([yshift=+1mm]box4.south) {};

\node[anchor=north] at (box4b.north) {%
    \begin{minipage}{18.8cm}
        \begin{itemize}[itemsep=4pt, topsep=4pt, left=5pt, label=\raisebox{0.35ex}{\scalebox{0.4}{$\blacksquare$}}]
            \item \textbf{TSA} with 933 periods achieves \textbf{zero error} in the OFV and aggregated decision variables
            \item \textbf{Lower bound on computational burden}: $51 \times 10^{-5}$~work units~$\rightarrow$~\textbf{369-fold} reduction versus full-scale model
        \end{itemize}
    \end{minipage}%
};

\end{tikzpicture}}
\caption{Conceptual four-step argument and accompanying example for exact TSA with ESSs.}
\label{fig:Concept}
\end{figure*}
{
\begin{subequations}
\label{eqn:model}
\begin{IEEEeqnarray}{l}
    \hspace{-2.5em}
    z_i = \min_{} \sum_{r\in\mathcal{R}^i}\!W_r \big(
          \sum_{g\in\mathcal{G}}  C_g              p_{r,g}              
        + \sum_{s\in\mathcal{S}}  C^\mathrm{d}_{s} p^\mathrm{d}_{r,s}    
        +                         C^\mathrm{ns}    p^\mathrm{ns}_{r}     
        \big)
    \label{eqn:of} \\
    \hspace{-2.5em}
    \sum_{g\in\mathcal{G}} p_{r,g} + \sum_{s\in\mathcal{S}} \big( p^{\mathrm{d}}_{r,s} - p^{\mathrm{c}}_{r,s} \big) + p^\mathrm{ns}_{r} = \widetilde{D}_{r}
    \quad (\smash{\mu}^\mathrm{bal}_{r}) \nonumber \\ \hspace{-2.5em}
    \quad \forall r\in\mathcal{R}^i
    \label{eqn:balance} \\
    \hspace{-2.5em}
    0 \leq p_{r,t} \leq \overline{P}_t\ 
    \quad (\smash{\underline{\lambda}}^\mathrm{t}_{r,t}, \smash{\overline{\lambda}}^\mathrm{t}_{r,t})
    \quad \forall r\in\mathcal{R}^i, t\in\mathcal{G}^\mathcal{T}
    \label{eqn:therm_bounds} \\
    \hspace{-2.5em}
    0 \leq p_{r,v} \leq \overline{P}_v \widetilde{F}_{r,v}
    \quad (\smash{\underline{\lambda}}^\mathrm{v}_{r,v}, \smash{\overline{\lambda}}^\mathrm{v}_{r,v})
    \quad \forall r\in\mathcal{R}^i, v\in\mathcal{G}^\mathcal{V}
    \label{eqn:res_bounds} \\
    \hspace{-2.5em}
    0 \leq p^{\mathrm{ns}}_{r} \leq \widetilde{D}_{r}
    \quad (\smash{\underline{\lambda}}^{\mathrm{ns}}_{r}, \smash{\overline{\lambda}}^{\mathrm{ns}}_{r})
    \quad \forall r\in\mathcal{R}^i
    \label{eqn:nsp_bounds} \\
    \hspace{-2.5em}
    e_{r,s} = \underline{E}_{s} + W_r \big( \eta^\mathrm{c}_{s}p^\mathrm{c}_{r,s} - p^\mathrm{d}_{r,s}/\eta^\mathrm{d}_{s} \big)
    \quad (\smash{\mu}^\mathrm{ini}_{r,s}) \nonumber \\ \hspace{-2.5em}
    \quad r\in\{1\},\forall s\in\mathcal{S}
    \label{eqn:soc_start} \\
    \hspace{-2.5em}
    e_{r,s} = e_{r-1,s} + W_r \big(\eta^\mathrm{c}_{s} p^\mathrm{c}_{r,s} - p^\mathrm{d}_{r,s}/\eta^\mathrm{d}_{s}\big)
    \quad (\smash{\mu}^\mathrm{intra}_{r,s}) \nonumber \\ \hspace{-2.5em}
    \quad \forall r\in\{2,\dots,R^i\}, s\in\mathcal{S}
    \label{eqn:soc_intra} \\
    \hspace{-2.5em}
    e_{r,s} = \underline{E}_{s}
    \quad (\smash{\mu}^\mathrm{fin}_{r,s})
    \quad r\in\{R^i\},\forall s\in\mathcal{S}
    \label{eqn:soc_fin} \\
    \hspace{-2.5em}
    \underline{E}_{s} \leq e_{r,s} \leq \overline{E}_{s}
    \quad (\smash{\underline{\lambda}}^\mathrm{soc}_{r,s}, \smash{\overline{\lambda}}^\mathrm{soc}_{r,s})
    \quad \forall r\in\mathcal{R}^i, s\in\mathcal{S}
    \label{eqn:soc_bounds} \\
    \hspace{-2.5em}
    0 \leq p^\mathrm{c}_{r,s} \leq \overline{P}^\mathrm{c}_{s}
    \quad (\smash{\underline{\lambda}}^\mathrm{c}_{r,s}, \smash{\overline{\lambda}}^\mathrm{c}_{r,s})
    \quad \forall r\in\mathcal{R}^i, s\in\mathcal{S}
    \label{eqn:sto_ch} \\
    \hspace{-2.5em}
    0 \leq p^\mathrm{d}_{r,s} \leq \overline{P}^\mathrm{d}_{s}
    \quad (\smash{\underline{\lambda}}^\mathrm{d}_{r,s}, \smash{\overline{\lambda}}^\mathrm{d}_{r,s})
    \quad \forall r\in\mathcal{R}^i, s\in\mathcal{S}
    \label{eqn:sto_dis} 
\end{IEEEeqnarray}
\end{subequations}
}

To capitalize on the potential for parallelization and improved computational performance of model~\eqref{eqn:model}, enabled by disaggregation and subsequent aggregation, the modeler must make two key decisions: (\textit{a}) how to disaggregate the full time horizon into submodels~$i$, that is, where to place the splits in the full time horizon, and (\textit{b}) how to aggregate periods $r$ within each submodel. 
However, time-coupling from ESS renders these decisions particularly challenging. In the following section, we establish the conditions for disaggregation and aggregation that inform these TSA decisions, thereby enabling the construction of an aggregated model that exactly preserves the optimal solution of its full-scale counterpart.

\section{Concept for Exact TSA with ESSs}
\label{sec:ExactTSA}
This section demonstrates that \textit{exact} TSA is achievable in the presence of ESS time-coupling constraints. This is accomplished by disaggregating the full-scale model into independent submodels according to each period's ACS, and by aggregating periods within each submodel using their corresponding dual information -- without incurring an error in model outputs -- as detailed in Sections~\ref{subsec:ExactDisaggregation} and~\ref {subsec:ExactAggregation}.

The goal of this section is purely theoretical: to demonstrate that \textit{exact} TSA is achievable under perfect information for an optimization model with ESS time-coupling constraints. To make this idea precise, we outline a conceptual argument -- structured in four steps -- that shows how an exact TSA could be constructed if the ACSs and duals were known in advance. This is not intended as an algorithm, but rather as a theoretical framework illustrating that exact TSA is attainable under idealized conditions. To make the exposition concrete, we illustrate these steps using a numerical example of model~\eqref{eqn:model}, detailed in Section~\ref{sec:NumRes}. Section~\ref{sec:ML} then shows how the same conceptual idea can be approximated in practice, when perfect information is unavailable.

In \textbf{Step 1}, we compute the full-scale model comprising 8736 periods, as a reference and to identify the ACSs and corresponding duals for each period. 

\textbf{Step 2} disaggregates the full-scale model into independent submodels using ACS information without introducing any error. This is achieved by placing the splits at periods where the state-of-charge lower bound~\eqref{eqn:soc_bounds} is active for two consecutive periods, $r-1$ and $r$, as thoroughly explained in Section~\ref{subsec:ExactDisaggregation}. Periods in which splits are placed represent independent single-period submodels referred to as unlinked submodels~$\mathcal{U}$, which comprise 1377 of the 8736 periods in the full-scale model. Submodels that span multiple periods are referred to as linked submodels~$\mathcal{L}$. For example, the first linked submodel spans 4 periods, the second one 5, and the longest 395. In total, 7359 out of 8736 periods are represented by linked submodels. Since all constructed submodels are independent, they can be solved in parallel, substantially improving computational performance while preserving the exact solution of their full-scale counterpart.

While disaggregating the model into submodels enables parallelization, \textbf{Step~3} outlines how ACS and dual information can be leveraged to inform aggregation within each submodel, further enhancing computational efficiency. By definition, unlinked submodels are not subject to time coupling. Consequently, as shown in~\cite{Wogrin2023}, they can be aggregated perfectly into periods that share the same ACS. Following this, the 1377 unlinked submodels obtained through Step~2 can be represented exactly by only 8 aggregated periods. In contrast, aggregation in linked submodels must preserve chronology. However, when adjacent periods share the same ACS and duals, they can be aggregated without loss of information, as detailed in Section~\ref{subsec:ExactAggregation}. For instance, the first linked submodel, which span 4 periods, can be exactly represented by 2 aggregated periods, while the longest one aggregates 395 periods into only 5. After aggregation, the longest submodel comprises 31 periods. In total, the 7359 periods represented by linked submodels can be compacted into just 925 periods.

In \textbf{Step 4}, all disaggregated and subsequently aggregated submodels are solved in parallel. With perfect information, we achieve a TSA with zero error in both the objective function value (OFV) and the aggregated decision variables. Due to parallelization, the lower bound on computational burden for solving the aggregated model (excluding parallelization overhead) results in a 369-fold reduction relative to the full-scale model. Fig.~\ref{fig:Concept} illustrates the conceptual four-step argument, including the accompanying numerical example representing one of the four co-scheduling problem instances considered, namely the BESS--Solar case, detailed in Section~\ref{sec:NumRes}. On average, the reduction in computational burden is 415-fold.

\subsection{Exact Disaggregation into Independent Submodels}
\label{subsec:ExactDisaggregation}

Previous work~\cite{Cardona2025} has shown that in models with ramping constraints (time-coupling inequality constraints), suitable periods for disaggregation occur where such constraints are inactive, because the system matrix could be decomposed exactly without incurring an error. In this paper, we extend these ideas to models with ESSs time-coupling equality constraints. 

In model~\eqref{eqn:model}, the state-of-charge constraint~\eqref{eqn:soc_intra} is the only time-coupling constraint. For illustration,~(2) presents the part of the system matrix that arises for two consecutive periods, i.e., $r-1,~r$, for constraint~\eqref{eqn:soc_intra}. 
Note that other constraints are omitted in the system matrix in this discussion because none of them are time-coupling. The vertical dashed line indicates the separation of variables of different periods. We observe that only state-of-charge variables $e_{r,s}$ are coupling consecutive periods, as indicated by red coefficients in the system matrix. By definition, an equality constraint is always active. However, if both $e_{r-1,s}$ and $e_{r,s}$ take the value zero, then the red coefficients in the system matrix can be ignored (because they are multiplied by zero) and the system matrix could be disaggregated into two separate submodels without incurring an error in model outputs.
\begin{figure}[t]
\centering
\input{matrix/matrix}
\label{fig:Matrix}
\end{figure}
Hence, the ESS being empty allows for disaggregation of the system matrix, and thus disaggregation of the full-scale model into submodels (i.e., Step 2). The ESS is considered empty, when the state-of-charge lower bound~\eqref{eqn:soc_bounds} is active for at least two consecutive periods, i.e., $e_{r-1,s}=e_{r,s}=0$. Periods in which the ESS is empty result in unlinked submodels.
To further reduce the size of the resulting submodels, future work will explore additional disaggregation criteria.

\subsection{Exact Temporal Aggregation within Submodels}
\label{subsec:ExactAggregation}
Previous research~\cite{Wogrin2023} has proven that in models without time-coupling constraints, periods sharing the same ACS can be aggregated exactly. This idea can be applied directly to the unlinked submodels in Step 3. It has also been extended to inequality ramping constraints~\cite{Cardona2024} when aggregating shorter time horizons that share the same ACSs in all periods. This paper is the first one to extend this idea to ESSs time-coupling equality constraints.
\begin{figure*}[t]
\centering
\input{acs/acs}
\caption{Disaggregation into submodels (top) and aggregation within submodels (bottom) via active constraint sets.}
\label{fig:States}
\end{figure*}

We have observed that when two consecutive periods share the same ACS and identical duals, aggregating these periods preserves both the ACS and the dual values -- hence, no information is lost through aggregation. In particular, the most relevant feature for the problem at hand appears to be the ACS characterization of the storage state (e.g., empty, charging, discharging, or idle) together with the identity of the marginal generator (i.e., the generator dispatched to supply an additional unit of power).
Intuitively, if the storage is charging over several consecutive periods and the marginal generator -- indicated by the dual of the balance constraint~\eqref{eqn:balance} -- remains unchanged, these periods can be collapsed into a single aggregated period without, on average, introducing error in the model outputs. A formal proof of this observation is left for future work; in this paper, we provide supporting evidence through the numerical examples in Section~\ref{sec:NumRes}.

To illustrate the concepts introduced in this section, Fig.~\ref{fig:States} shows the disaggregation of the full time horizon into unlinked and linked submodels (Step~2). Periods when the ESS is empty yield unlinked submodels, as indicated by ACS1 in the top panel. The bottom panel illustrates Step~3, where periods sharing the same ACS and marginal generator -- indicated by the same color -- are aggregated. This aggregation preserves temporal chronology across ACSs while merging periods within each ACS.

\section{ML Framework for TSA with ESSs}
\label{sec:ML}
The previous section outlined how the co-scheduling model can be disaggregated into independent, parallelizable submodels via ACS identification, followed by temporal aggregation to enhance computational efficiency. The identification of ACSs depends on a priori knowledge of the model’s optimal solution space. To overcome this limitation and enable a practical method to disaggregation and aggregation, we build on the theoretical insights from the preceding analysis to develop an ML framework that aim to achieve exact disaggregation (Section~\ref{sec:ML-Disaggregation}) and aggregation (Section~\ref{sec:ML-Aggregation}) of the co-scheduling model.

\subsection{Practical Disaggregation into Independent Submodels}
\label{sec:ML-Disaggregation}
We now describe the ML-based disaggregation of the full time horizon. Specifically, full-scale model solutions over a 3-year horizon, together with the corresponding input time series, are used for feature identification prior to training a classifier for a given power system configuration, which can then be applied to different realizations of the input data.
Each classifier is tested using entirely new input data comprising a 1\nobreakdash-year horizon and evaluated against the solution of the corresponding full-scale model.

\subsubsection{Feature Identification}
\label{sec:Features}
First, we employ logistic regression to identify the most informative features from the candidate set described in Table~\ref{tab:CandFeatures}, which also includes lead and lagged versions up to $\pm 10$ periods, a design choice that demonstrated good performance across the power system configuration considered in our numerical experiments. Feature importance is quantified from the regression coefficients, and only those contributing more than 1\% (following the elbow criterion) are retained for classifier training.
\begin{table}[t]
\renewcommand{\arraystretch}{1}
\centering
\caption{Candidate Features.}
\label{tab:CandFeatures}
\begin{tabular}{l l}
\hline
Feature & Description \\ \hline
$\textit{D}$                                & Demand \\
$\textit{F}$                                & Capacity factor \\
$\overline{\textit{VRE}}$                   & Capacity factor times VRE capacity \\
\textit{$\overline{\textit{VRE}}$-D}        & Ratio of $\overline{\textit{VRE}}$ and $D$ \\
$\textit{Crit}$                             & Ratio of $D$ and total generation capacity ($\overline{\textit{VRE}}$ and thermal) \\
\hline
\end{tabular}
\end{table}
\begin{table}[t]
\renewcommand{\arraystretch}{1}
\centering
\caption{Classifier Performance -- Validation and Test Data.}
\label{tab:CLPerf}
\begin{tabular}{l l l l}
\hline
Dataset & Truth     & \multicolumn{2}{l}{Prediction} \\ \cline{3-4}
(Accuracy, Balanced Accuracy)    &       & Empty & Not Empty \\ \hline
\multirow{2}{*}{\makecell[l]{Validation \\ (92\%, 84\%)}} 
        & Empty   & 690  & 292  \\
        & Not Empty & 113  & 4149 \\ \hline
\multirow{2}{*}{\makecell[l]{Test \\ (86\%, 69\%)}} 
        & Empty   & 642 & 873 \\
        & Not Empty & 339 & 6882 \\ \hline
\end{tabular}
\end{table}

\subsubsection{Classification}
\label{sec:Class}
Second, for each power system configuration we train a binary classifier -- using the relevant features determined from logistic regression -- to identify the periods when the ESS is empty, i.e., when the state-of-charge lower bound~\eqref{eqn:soc_bounds} is active for at least two consecutive periods, as detailed in Section~\ref{subsec:ExactDisaggregation}. The classifier uses a random forest as this has shown excellent performance in previous work~\cite{Cardona2025}. For training the classifier we use an 80-20 split of the 3-year time horizon data, with 80\% used for model fitting and 20\% for hyperparameter tuning via k-fold cross-validation (i.e., validation).
The set-up of the ML-based disaggregation, comprising feature identification and classifier training, requires approximately 2~hours per system configuration and is performed offline. The classifier is then tested on 1\nobreakdash-year out-of-sample datasets.
As shown in Table~\ref{tab:CLPerf}, the classifier maintains a desirable level of accuracy on both the validation and test datasets (BESS–Solar case).
Misclassification has two effects: false positives (i.e., periods where the ESS is not empty but classified as empty) cause minor deviations in the disaggregated model outcomes, whereas false negatives (i.e., periods where the ESS is empty but classified as not empty) lead to longer submodels, leaving disaggregation potential untapped.

For the ML-based disaggregation, we used full-scale model solutions to obtain the ACS for feature identification and classifier training and validation. Nevertheless, the approach remains highly relevant for practical applications, as it could, for example, rely on partial model runs, thereby alleviating the need to solve the full-scale problem -- a direction we plan to explore in future research. Moreover, once trained for a given system configuration, the ML-based disaggregation could be applied in Monte Carlo–type analyses, in which the co-scheduling problem is repeatedly solved for many realizations of uncertain parameters, such as renewable capacity factors or demand.
\begin{table*}[t]
\renewcommand{\arraystretch}{1.2}
\centering
\setlength{\tabcolsep}{3pt} 
\caption{ML-based Disaggregation Results Relative to Full-Scale Model.}
\label{tab:DisaggStatistics}
\begin{tabular}{lllllllllllllll}
\hline
Case        & \multicolumn{3}{l}{Submodel statistics}                             & & \multicolumn{6}{l}{Relative error in \%}                      & & \multicolumn{2}{l}{Comp. burden in work units$\times 10^{-3}$}  & Speed-up in p.u.  \\\cline{2-4}\cline{6-11} \cline{13-14}
~           & $|\mathcal{U}|$ & $|\mathcal{L}|$ & $\max(\mathcal{R}^\mathcal{L})$ & & OFV  & VRE   & Thermal & NSP     & Charging   & Discharging   & &  Full-scale & Largest submodel             &                         \\ \hline
BESS--Solar & 993             &  239            & 405                             & & 0.27 &  0.00 &  0.00   &    7.65 &   0.40     &   0.40        & &  187.21     &  9.79                        & 19.12                   \\
BESS--Wind  & 797             &  196            & 219                             & & 1.18 & -0.03 &  0.03   &   24.86 &   0.39     &   0.39        & &  198.28     &  7.45                        & 26.61                   \\
PHS--Solar  & 120             &  71             & 3350                            & & 7.74 & -0.01 & -0.08   & 7696.83 &   1.39     &   1.39        & &  388.50     & 64.35                        &  6.04                   \\
PHS--Wind   & 804             &  183            & 335                             & & 2.55 & -0.97 &  1.32   &  734.52 & -13.87     & -13.87        & &  251.48     &  9.91                        & 25.37                   \\
\hline
\end{tabular}
\end{table*}
\begin{figure*}[t]
\centering
\definecolor{ErrObj}{HTML}{C62324}
\definecolor{ErrVRE}{HTML}{6A9E1F}
\definecolor{ErrThermal}{HTML}{E87D37}
\definecolor{ErrNSP}{HTML}{052F61}
\definecolor{ErrCharging}{HTML}{14967C}
\definecolor{ErrDischarging}{HTML}{A50E82}

\begin{tikzpicture}
\begin{groupplot}[
    name=mygroupplot, 
    group style={
        group size=2 by 2,
        horizontal sep=1.5cm,
        vertical sep=1.2cm,
        xlabels at=edge bottom,
        ylabels at=edge left
    },
    width=0.85\columnwidth,
    height=3.5cm,
    scale only axis,
    ymin=-50, ymax=10,
    xmin=30, xmax=100,
    xtick={10,20,30,40,50,60,70,80,90,100},
    tick style={major tick length=3pt, semithick},
    ytick={-50,-40,-30,-20,-10,0,10},
    xlabel={Clusters},
    ylabel={Relative error in \%},
    grid=none
]

\nextgroupplot[title={\textbf{BESS--Solar}}, title style={align=center, yshift=-2mm}]
\addplot[color=ErrVRE, mark=*, thick] table[x=CL,y=VRE,col sep=comma]{results/RelError_Hier_NetD_BESS_Solar.csv};
\addplot[color=ErrThermal, mark=diamond*, thick] table[x=CL,y=Thermal,col sep=comma]{results/RelError_Hier_NetD_BESS_Solar.csv};
\addplot[color=ErrNSP, mark=square*, thick] table[x=CL,y=NSP,col sep=comma]{results/RelError_Hier_NetD_BESS_Solar.csv};
\addplot[color=ErrCharging, mark=pentagon*, thick] table[x=CL,y=Ch,col sep=comma]{results/RelError_Hier_NetD_BESS_Solar.csv};
\addplot[color=ErrDischarging, mark=triangle*, thick] table[x=CL,y=Dis,col sep=comma]{results/RelError_Hier_NetD_BESS_Solar.csv};
\addplot[color=ErrObj, mark=asterisk, thick] table[x=CL,y=OFV,col sep=comma]{results/RelError_Hier_NetD_BESS_Solar.csv};

\nextgroupplot[title={\textbf{PHS--Solar}}, title style={align=center, yshift=-2mm}]
\addplot[color=ErrVRE, mark=*, thick] table[x=CL,y=VRE,col sep=comma]{results/RelError_Hier_NetD_PHS_Solar.csv};
\addplot[color=ErrThermal, mark=diamond*, thick] table[x=CL,y=Thermal,col sep=comma]{results/RelError_Hier_NetD_PHS_Solar.csv};
\addplot[color=ErrNSP, mark=square*, thick] table[x=CL,y=NSP,col sep=comma]{results/RelError_Hier_NetD_PHS_Solar.csv};
\addplot[color=ErrCharging, mark=pentagon*, thick] table[x=CL,y=Ch,col sep=comma]{results/RelError_Hier_NetD_PHS_Solar.csv};
\addplot[color=ErrDischarging, mark=triangle*, thick] table[x=CL,y=Dis,col sep=comma]{results/RelError_Hier_NetD_PHS_Solar.csv};
\addplot[color=ErrObj, mark=asterisk, thick] table[x=CL,y=OFV,col sep=comma]{results/RelError_Hier_NetD_PHS_Solar.csv};

\nextgroupplot[title={\textbf{BESS--Wind}}, title style={align=center, yshift=-2mm}]
\addplot[color=ErrVRE, mark=*, thick] table[x=CL,y=VRE,col sep=comma]{results/RelError_Hier_NetD_BESS_Wind.csv};
\addplot[color=ErrThermal, mark=diamond*, thick] table[x=CL,y=Thermal,col sep=comma]{results/RelError_Hier_NetD_BESS_Wind.csv};
\addplot[color=ErrNSP, mark=square*, thick] table[x=CL,y=NSP,col sep=comma]{results/RelError_Hier_NetD_BESS_Wind.csv};
\addplot[color=ErrCharging, mark=pentagon*, thick] table[x=CL,y=Ch,col sep=comma]{results/RelError_Hier_NetD_BESS_Wind.csv};
\addplot[color=ErrDischarging, mark=triangle*, thick] table[x=CL,y=Dis,col sep=comma]{results/RelError_Hier_NetD_BESS_Wind.csv};
\addplot[color=ErrObj, mark=asterisk, thick] table[x=CL,y=OFV,col sep=comma]{results/RelError_Hier_NetD_BESS_Wind.csv};

\nextgroupplot[name=lastsubplot, title={\textbf{PHS--Wind}}, title style={align=center, yshift=-2mm},
    legend to name=commonlegend,    
    legend columns=6,
    legend style={
        draw=none,
        column sep=0.2cm,          
        nodes={scale=1, transform shape, inner xsep=0mm} 
    }
]
\addplot[color=ErrVRE, mark=*, thick] table[x=CL,y=VRE,col sep=comma]{results/RelError_Hier_NetD_PHS_Wind.csv};
\addlegendentry{VRE}

\addplot[color=ErrThermal, mark=diamond*, thick] table[x=CL,y=Thermal,col sep=comma]{results/RelError_Hier_NetD_PHS_Wind.csv};
\addlegendentry{Thermal}

\addplot[color=ErrNSP, mark=square*, thick] table[x=CL,y=NSP,col sep=comma]{results/RelError_Hier_NetD_PHS_Wind.csv};
\addlegendentry{NSP}

\addplot[color=ErrCharging, mark=pentagon*, thick] table[x=CL,y=Ch,col sep=comma]{results/RelError_Hier_NetD_PHS_Wind.csv};
\addlegendentry{Charging}

\addplot[color=ErrDischarging, mark=triangle*, thick] table[x=CL,y=Dis,col sep=comma]{results/RelError_Hier_NetD_PHS_Wind.csv};
\addlegendentry{Discharging}

\addplot[color=ErrObj, mark=asterisk, thick] table[x=CL,y=OFV,col sep=comma]{results/RelError_Hier_NetD_PHS_Wind.csv};
\addlegendentry{OFV}

\end{groupplot}

\node at ($(lastsubplot.south)+(-4.5cm,9.2cm)$) 
    {\pgfplotslegendfromname{commonlegend}};

\begin{axis}[
    at={(group c1r1.north east)},   
    anchor=north east,
    xshift=-0.3cm,
    yshift=-1.1cm,               
    width=4.6cm,
    height=3.5cm,
    xmin=30,
    xmax=70,                      
    ymin=-4,
    ymax=2,
    grid=none,
    xtick={10,20,30,40,50,60,70},
    ytick={-8,-6,-4,-3,-2,-1,0,1,2},
    tick label style={font=\tiny},
    label style={font=\scriptsize},
    tick style={major tick length=2pt, semithick}
]

\addplot[color=ErrVRE, mark=*, mark size=1pt, thick] 
    table[x=CL,y=VRE,col sep=comma] {results/RelError_Hier_NetD_BESS_Solar.csv};

\addplot[color=ErrThermal, mark=diamond*, mark size=1pt, thick] 
    table[x=CL,y=Thermal,col sep=comma] {results/RelError_Hier_NetD_BESS_Solar.csv};

\addplot[color=ErrNSP, mark=square*, mark size=1pt, thick] 
    table[x=CL,y=NSP,col sep=comma] {results/RelError_Hier_NetD_BESS_Solar.csv};

\addplot[color=ErrCharging, mark=pentagon*, mark size=1pt, thick] 
    table[x=CL,y=Ch,col sep=comma] {results/RelError_Hier_NetD_BESS_Solar.csv};

\addplot[color=ErrDischarging, mark=triangle*, mark size=1pt, thick] 
    table[x=CL,y=Dis,col sep=comma] {results/RelError_Hier_NetD_BESS_Solar.csv};

\addplot[color=ErrObj, mark=asterisk, mark size=1.5pt, thick] 
    table[x=CL,y=OFV,col sep=comma] {results/RelError_Hier_NetD_BESS_Solar.csv};
\end{axis}

\begin{axis}[
    at={(group c1r2.north east)},   
    anchor=north east,
    xshift=-0.3cm,
    yshift=-1.1cm,               
    width=4.6cm,
    height=3.5cm,
    xmin=30,
    xmax=70,                      
    ymin=-4,
    ymax=2,
    grid=none,
    xtick={10,20,30,40,50,60,70},
    ytick={-8,-6,-4,-3,-2,-1,0,1,2},
    tick label style={font=\tiny},
    label style={font=\scriptsize},
    tick style={major tick length=2pt, semithick}
]

\addplot[color=ErrVRE, mark=*, mark size=1pt, thick] 
    table[x=CL,y=VRE,col sep=comma] {results/RelError_Hier_NetD_BESS_Wind.csv};

\addplot[color=ErrThermal, mark=diamond*, mark size=1pt, thick] 
    table[x=CL,y=Thermal,col sep=comma] {results/RelError_Hier_NetD_BESS_Wind.csv};

\addplot[color=ErrNSP, mark=square*, mark size=1pt, thick] 
    table[x=CL,y=NSP,col sep=comma] {results/RelError_Hier_NetD_BESS_Wind.csv};

\addplot[color=ErrCharging, mark=pentagon*, mark size=1pt, thick] 
    table[x=CL,y=Ch,col sep=comma] {results/RelError_Hier_NetD_BESS_Wind.csv};

\addplot[color=ErrDischarging, mark=triangle*, mark size=1pt, thick] 
    table[x=CL,y=Dis,col sep=comma] {results/RelError_Hier_NetD_BESS_Wind.csv};

\addplot[color=ErrObj, mark=asterisk, mark size=1.5pt, thick] 
    table[x=CL,y=OFV,col sep=comma] {results/RelError_Hier_NetD_BESS_Wind.csv};
\end{axis}

\begin{axis}[
    at={(lastsubplot.north east)},   
    anchor=north east,
    xshift=-0.3cm,
    yshift=-1.1cm,               
    width=4.6cm,
    height=3.5cm,
    xmin=30,
    xmax=70,                      
    ymin=-2,
    ymax=4,
    grid=none,
    xtick={10,20,30,40,50,60,70},
    ytick={-8,-6,-4,-3,-2,-1,0,1,2,3,4},
    tick label style={font=\tiny},
    label style={font=\scriptsize},
    tick style={major tick length=2pt, semithick}
]

\addplot[color=ErrVRE, mark=*, mark size=1pt, thick] 
    table[x=CL,y=VRE,col sep=comma] {results/RelError_Hier_NetD_PHS_Wind.csv};

\addplot[color=ErrThermal, mark=diamond*, mark size=1pt, thick] 
    table[x=CL,y=Thermal,col sep=comma] {results/RelError_Hier_NetD_PHS_Wind.csv};

\addplot[color=ErrNSP, mark=square*, mark size=1pt, thick] 
    table[x=CL,y=NSP,col sep=comma] {results/RelError_Hier_NetD_PHS_Wind.csv};

\addplot[color=ErrCharging, mark=pentagon*, mark size=1pt, thick] 
    table[x=CL,y=Ch,col sep=comma] {results/RelError_Hier_NetD_PHS_Wind.csv};

\addplot[color=ErrDischarging, mark=triangle*, mark size=1pt, thick] 
    table[x=CL,y=Dis,col sep=comma] {results/RelError_Hier_NetD_PHS_Wind.csv};

\addplot[color=ErrObj, mark=asterisk, mark size=1.5pt, thick] 
    table[x=CL,y=OFV,col sep=comma] {results/RelError_Hier_NetD_PHS_Wind.csv};
\end{axis}

\end{tikzpicture}
\caption{Solution error of ML-based aggregation for the considered system configurations, relative to the disaggregated model.}
\label{fig:RelError}
\end{figure*}

\subsection{Practical Temporal Aggregation within Submodels}
\label{sec:ML-Aggregation}
In~\cite{santosuosso2025clusteringforestablishingperformance}, we demonstrated that an appropriately constructed aggregated model consistently provides a valid lower bound on the optimal OFV of its full-scale counterpart. In Section~\ref{sec:ExactTSA}, we extend this result by showing that the lower bound converges to the optimal OFV of the disaggregated model.
Building on this, we aggregate periods within linked submodels -- identified in Section~\ref{sec:ML-Disaggregation} -- while preserving chronology, as this was a necessary condition for perfect aggregation. For simplicity, the ML-based aggregation follows a hierarchical clustering method where adjacent periods are clustered based on the similarity of the net demand~$D_r^{\mathrm{net}}$, i.e., $D_r^{\mathrm{net}}=\max\{0,\widetilde{D}_r-\sum_{v\in\mathcal{G}^{\mathcal{V}}}\overline{P}_v\widetilde{F}_{r,v}\}~\forall~r\in\mathcal{R}^i$, until the number of periods is reduced to the desired number of clusters. 
The intuition behind this approach is that periods with comparable net demand are likely to correspond to similar system operating conditions -- and are therefore characterized by the same ACS and duals. Aggregating such periods is expected to approximate the behavior of perfect aggregation in practice. In future work, we will refine the ML-based aggregation using a nested ML framework informed by ACS and dual information.

\section{Numerical Results}
\label{sec:NumRes}
This section presents numerical results for TSA of four instances of the co-scheduling model~\eqref{eqn:model}, following the ML framework described in Section~\ref{sec:ML}.

Each of the four cases represents a 1-year time horizon of 8736 hours and includes one ESS -- either a battery energy storage system (BESS) or a pumped hydro storage (PHS) unit -- together with a VRE source -- either solar or wind. In addition, all cases include a thermal unit and allow non-supplied power (NSP) to ensure that demand is always met.
The BESS has a capacity of $\overline{E}_s = 400\,\text{MWh}$, charging and discharging power of $\smash{\overline{P}}^\mathrm{c}_s = \smash{\overline{P}}^\mathrm{d}_s = 100\,\text{MW}$, with corresponding efficiencies of $\eta^{\mathrm{c}}_s = \eta^{\mathrm{d}}_s = 0.92$. The variable cost of discharging is $C^{\mathrm{d}}_s = \text{\texteuro}\,1.5\,\text{MWh}^{-1}$.
For the PHS unit, the parameters are $\overline{E}_s = 1600\,\text{MWh}$, $\smash{\overline{P}}^\mathrm{c}_s = \smash{\overline{P}}^\mathrm{d}_s = 100\,\text{MW}$, and $\eta^{\mathrm{c}}_s = \eta^{\mathrm{d}}_s = 0.9$, with $C^{\mathrm{d}}_s = \text{\texteuro}\,0.5\,\text{MWh}^{-1}$.
The maximum generation of solar and wind units is $\overline{P}_v = 1000\,\text{MW}$, and that of the thermal units is $\overline{P}_t = 480\,\text{MW}$, with variable operating costs $C_g$ of $\text{\texteuro}\,1,~\text{\texteuro}\,2.5,~\text{and}~\text{\texteuro}\,60\,\text{MWh}^{-1}$, respectively. Finally, the cost of non-supplied energy is assumed as $C^{\mathrm{ns}} = \text{\texteuro}\,5000\,\text{MWh}^{-1}$. 
Time series comprise the downscaled Austrian demand, totaling 3363 GWh, and capacity factors from Renewables Ninja, which were excluded from classifier training, validation, and testing.

First, we use ML-based disaggregation to predict where to disaggregate the full time horizon into submodels based solely on input data. Table~\ref{tab:DisaggStatistics} presents the corresponding results, including the number of unlinked and linked submodels, the maximum length of linked submodels, the resulting output errors, the computational burden of solving the full-scale model compared with the largest submodel, as well as the resulting computational speed-up (excluding parallelization overhead).
The ML-based disaggregation delivers accurate results for the BESS cases, with OFV errors below 1.5\%.
These deviations primarily stem from an imperfect representation of NSP, which corresponds to a critical system state occurring in only 21 and 24 periods, respectively.
Specifically, the total NSP in the full-scale solutions amounts to 993~MWh and 830~MWh for the BESS--Solar and BESS--Wind cases, respectively, compared with 1069~MWh and 1036~MWh in the corresponding disaggregated models.
In the PHS cases, total NSP in the full-scale models is marginal, amounting to only 27~MWh and 31~MWh and occurring in just four periods for the PHS--Solar and PHS--Wind cases, respectively. In contrast to the BESS cases, NSP deviations under the disaggregated models are considerably larger, yielding 2116~MWh and 257~MWh, which in turn lead to differences in the OFV.
This behavior is not unexpected, as classification performance typically deteriorates for highly infrequent events, particularly when the classifier is not explicitly designed and trained to detect them.
Nevertheless, the proposed ML-based disaggregation substantially reduces the computational burden by enabling parallelization of submodels. This results in a 19-fold average speed-up across the studied cases when comparing the full-scale model with the largest submodel.

Second, after the ML-based disaggregation, periods within the linked submodels are further aggregated by employing ML-based aggregation.
Fig.~\ref{fig:RelError} reports the approximation error incurred by the ML-based aggregation relative to the disaggregated models of the previous step, as the number of clusters increases. Except PHS--Solar, the relative errors in OFV, VRE and thermal generation, charging, and discharging remain within $\pm5\%$ even with only 30 clusters, while NSP stabilizes at around 70 clusters.

Finally, we compare the performance of the proposed ML framework, which combines ML-based disaggregation and aggregation,
with an a priori aggregation-only baseline in which TSA is based solely on hierarchical clustering without prior disaggregation. Table~\ref{tab:PerformanceComp} reports the OFV error and computational speed-up for the proposed ML framework and the aggregation-only baseline relative to the full-scale model for the co-scheduling problem instances, considering 30 and 100 clusters. 
In all cases, our ML framework substantially outperforms the aggregation-only baseline in OFV accuracy, at the cost of a slight increase in the computational complexity of the temporally aggregated model, while still achieving average speed-ups of 486-fold and 119-fold for 30 and 100 clusters, respectively, thereby remaining significantly more efficient than solving the full-scale model directly.

\begin{table*}[!t]
\renewcommand{\arraystretch}{1.2}
\centering
\setlength{\tabcolsep}{3pt} 
\caption{Performance Comparison of Proposed ML Framework and Aggregation-only Baseline Relative to Full-Scale Model.}
\label{tab:PerformanceComp}
\begin{tabular}{lllllllllllll}
\hline
                &           & \multicolumn{5}{l}{OFV error in M€ (\%)}                                          & & \multicolumn{5}{l}{Speed-up in p.u.}                                        \\ \cline{3-7} \cline{9-13} 
                &           & \multicolumn{2}{l}{ML framework}     & & \multicolumn{2}{l}{Aggregation-only}     & & \multicolumn{2}{l}{ML framework}     & & \multicolumn{2}{l}{Aggregation-only}   \\ \cline{3-4} \cline{6-7} \cline{9-10} \cline{12-13}
Case            & Clusters  &    30             & 100              & &  30               & 100                  & & 30      & 100                        & & 30      & 100                      \\ \hline
BESS--Solar     &           &    -2.01 (-1.44)  & 0.35 (0.25)      & &  -15.03 (-10.74)  & -14.70 (-10.50)      & & 343.92  & 91.07                      & & 423.69  & 127.54                   \\
BESS--Wind      &           &     0.81  (0.91)  & 1.05 (1.18)      & &  -22.90 (-25.73)  & -16.02 (-18.01)      & & 356.35  & 84.40                      & & 485.94  & 131.96                   \\
PHS--Solar      &           &     4.25  (3.18)  & 6.03 (4.51)      & &  -8.80   (-6.58)  & -8.62   (-6.44)      & & 767.07  & 181.59                     & & 825.42  & 248.27                   \\
PHS--Wind       &           &     1.57  (1.92)  & 2.00 (2.45)      & &  -15.86 (-19.43)  & -9.89  (-12.12)      & & 476.82  & 118.72                     & & 587.61  & 159.41                   \\
\hline
\end{tabular}
\end{table*}

\section{Conclusions}
\label{sec:Conclusions}
This paper showed that \textit{exact} temporal aggregation of optimization models with energy storage time-coupling constraints is possible, representing a key methodological advancement for the field. We presented a conceptual argument showing that the time horizon of a full-scale model can be disaggregated into independent, parallelizable submodels, with further aggregation applied within these submodels. Numerical experiments validated the theoretical results: under perfect information, ACSs and dual variables obtained from the full-scale model solution enable disaggregation and aggregation while exactly preserving the objective function value and decision variables, resulting in a theoretical 415-fold average reduction in computational burden in our case study. To bridge our theoretical analysis and practical application, an ML framework combining disaggregation and aggregation was developed, achieving average speed-ups of up to 486-fold for the considered co-scheduling problem instances.
A comparison with an a priori aggregation-only baseline shows that our ACS-informed ML framework yields aggregated models that more accurately capture the optimal OFV of their full-scale counterparts, while maintaining a favorable trade-off between solution accuracy and complexity.

While these results are promising, they open several directions for future research.
For example, extending the proposed ML framework to estimate the observed ACSs and duals at each period represents a promising direction for further improving the accuracy of ML-based aggregation. Moreover, extending the proposed solution framework to incorporate network constraints, ramping constraints, multiple ESSs, and long-duration ESSs constitutes another direction. In this regard, ACS-informed aggregation has already been successfully utilized for optimal power flow models~\cite{Stoeckl2025}, highlighting its applicability to a broader class of power system optimization problems than those considered here.

\bibliographystyle{IEEEtran}
\bibliography{main}

\end{document}